\documentclass[12pt,reqno,twoside]{amsart}
\usepackage{amscd,amsthm}
\IfFileExists{euler.sty}{\IfFileExists{palatino.sty}{\usepackage{palatino}}{}}{}

\makeatletter
\@ifpackageloaded{palatino}{\usepackage[mathcal]{euler}}{}
\makeatother

\AtEndDocument{
  \vfill\eject
  \ifdim\paperheight=11in
   \message{RECOMMENDATION: If your paper is A4, insert the "a4paper" 
option 
   in the "documentclass" command.}
   \fi
}

\numberwithin{equation}{section}

\newtheorem{thm}{Theorem}[section]
\newtheorem{dthm}[thm]{Definition/Theorem}
\newtheorem{lem}[thm]{Lemma}
\newtheorem{cor}[thm]{Corollary}

\theoremstyle{definition}
\newtheorem*{definition}{Definition}

\def\beq#1\eeq{\begin{equation}#1\end{equation}}
\makeatletter
\@ifpackageloaded{euler}{
 \newcommand{\fchoice}[2]{#1}
 \newcommand{\onto}{\to\mkern-14mu\to}
 \def\rightarrowfill@#1{\m@th\setboxz@h{$#1\relbar$}\ht\z@\z@
   $#1\copy\z@\mkern-6mu\cleaders
   \hbox{$#1\mkern-2mu\box\z@\mkern-2mu$}\hfill
   \mkern-6mu\mathord\rightarrow$}
 \def\leftarrowfill@#1{\m@th\setboxz@h{$#1\relbar$}\ht\z@\z@
   $#1\mathord\leftarrow\mkern-6mu\cleaders
   \hbox{$#1\mkern-2mu\copy\z@\mkern-2mu$}\hfill
   \mkern-6mu\box\z@$}
 \def\B@R#1#2{\raisebox{-.07ex}{$#1#2$}\mkern-6mu}
 \renewcommand{\hbar}{{\mspace{1mu}\mathpalette\B@R{\mathchar'26}h}}
}{
 \AtEndDocument{\IfFileExists{palatino.sty}{\message{RECOMMENDATION: 
 Get the "euler" package from tug.ctan.org.}}{\message{RECOMMENDATION: 
 Get the "palatino" and "euler" packages from tug.ctan.org.  }}}
 \newcommand{\fchoice}[2]{#2}
 \DeclareMathSymbol{\onto}{\mathrel}{AMSa}{"10}
}
\makeatother
\DeclareMathSymbol{\smallsetminus}{\mathbin}{AMSb}{"72}
\newcommand{\A}{\mathcal{A}}
\renewcommand{\AA}{\mathbb{A}}
\newcommand{\VV}{\mathbb{V}}
\newcommand{\WW}{\mathbb{W}}
\newcommand{\C}{\mathcal{C}}
\newcommand{\M}{\mathcal{M}}
\newcommand{\co}{\mathbb{C}}
\newcommand{\I}{\mathcal{I}}
\newcommand{\N}{\mathbb{N}}
\newcommand{\cs}{C${}^*$} 
\newcommand{\Z}{\mathbb{Z}} 
\newcommand{\into}{\hookrightarrow}
\newcommand{\isom}{\mathrel{\widetilde\longrightarrow}}
\newcommand{\Po}{\mathcal{P}}
\newcommand{\Gh}{\Gamma_{\!\mathrm{hol}}}
\newcommand{\HN}{\mathcal{H}_N} 
\newcommand{\LN}{L_N} 
\newcommand{\jb}{{\bar{\jmath}}} 
\newcommand{\ib}{{\bar{\imath}}} 
\DeclareMathOperator{\End}{End} 
\DeclareMathOperator{\Hom}{Hom}
\DeclareMathOperator{\RealPart}{Re}
 \renewcommand{\Re}{\RealPart}
\DeclareMathOperator{\Image}{Im}
 \renewcommand{\Im}{\Image}
\DeclareMathOperator{\id}{id}
\DeclareMathOperator{\rk}{rk}
\DeclareMathOperator{\td}{td}
\DeclareMathOperator{\ch}{ch}
\DeclareMathOperator{\tr}{tr}
\newcommand{\ntr}{\widetilde{\smash\tr\vphantom{\raise.5pt\hbox{r}}}}
\DeclareMathOperator{\Tr}{Tr}

\DeclareMathOperator{\coker}{coker}
\DeclareMathOperator{\spec}{Spec}
\DeclareMathOperator{\vol}{vol}
\newcommand{\Kahler}{K\"ahler} 
\newcommand{\norm}[1]{\lVert#1\rVert}
\newcommand{\Norm}[1]{\left\|#1\right\|}
\newcommand{\lef}{E_N}
\newcommand{\func}{\mathcal{Q}}
\newcommand{\gq}{Q}
\newcommand{\intM}{\int_\M\mathchoice{\mskip-2\thinmuskip}%
                   {\mskip-\thinmuskip}{}{}}
\newcommand{\TM}{T\fchoice{}{\!}\M}

\title{Geometric Quantization of Vector Bundles}
\author{Eli Hawkins}

\begin{document}
\begin{flushright}
\vspace*{-0.4in}
\begin{tabular}{l}
\textsf{\small CGPG-98/8-2}\\
\textsf{\small math.QA/9808116}\\
\end{tabular}
\vspace{0.25in}
\end{flushright}

\maketitle
\begin{center}
\vspace{-4ex}
\small\emph{\small Center for Gravitational Physics and Geometry}\\
\emph{\small The Pennsylvania State University,
University Park, PA 16802}\\
{\small E-mail: mrmuon@gravity.phys.psu.edu}\\
\end{center}

\begin{abstract}
I repeat my definition for quantization of a vector bundle.  For the 
cases of the Toeplitz and geometric quantizations of a compact \Kahler\ 
manifold, I give a construction for quantizing any smooth vector 
bundle which depends functorially on a choice of connection on the 
bundle.
\end{abstract}

\section*{Introduction}
Traditionally, ``quantization'' has meant some sort of process that, 
given a classical, symplectic phase space, produces a noncommutative 
algebra of quantum observables.  The concept of noncommutative 
geometry (see \cite{con1}) suggests that such a noncommutative algebra 
can be thought of as the algebra of functions on a ``noncommutative 
space''\!, so perhaps quantization could be made into a way of 
constructing a noncommutative geometry from a classical geometry.

However, as it stands, quantization is \emph{only} a procedure for 
constructing an algebra.  Since the algebra of continuous (or smooth) 
functions contains only the information of the point-set (or 
differential) topology of a space, this is merely the quantization of 
topology.  It would be desirable to extend quantization to a theory of 
quantization of \emph{geometry}.

Beyond topology, vector bundles are arguably the second most 
fundamental structure in geometry, so a plausible first step towards a 
theory of quantizing geometry would be a theory of quantizing vector 
bundles.  I began constructing such a theory in \cite{haw1} by giving 
a definition of vector bundle quantization and a procedure for 
quantizing the equivariant vector bundles over coadjoint orbits of 
compact, semi-simple Lie groups.

I continue the story in this paper by giving a procedure for 
quantizing arbitrary smooth vector bundles over compact \Kahler\ 
manifolds.  The construction depends only on the structures used to 
quantize the manifold, the vector bundle itself, and a connection on 
the vector bundle.

\section{Generalities}
\label{general}
Most of the symbols defined in this section will be defined again 
through constructions in Sec's \ref{Toplitz} and \ref{convergence}.  
The theorems of Sec.~\ref{convergence} will show that these 
constructions actually do satisfy the original definitions.

Recall that a continuous field of \cs-algebras (see \cite{dix1,k-w1}) 
is the natural notion of a bundle of \cs-algebras.  The fibers are 
all \cs-algebras, the space of (continuous) sections is a \cs-algebra, 
and for each point of the base space there is an evaluation map, a 
$*$-homomorphism of the algebra of sections onto the fiber algebra.

By the most general definition (see \cite{lan1,lan2,rie1}), a strict 
deformation quantization of a (Poisson) manifold $\M$ consists of a 
continuous field of \cs-algebras, $\A_{\hat\I}$, and a (total) 
quantization map.  Conventionally, the base space $\hat\I$ of the 
continuous field is the set of possible values of $\hbar$; more 
generally, it is just some set containing an accumulation point 
$\infty\in\hat\I$ which plays the role of $\hbar=0$.  The fiber of the 
continuous field at this ``classical limit'' point is $\C(\M)$, the 
\cs-algebra of continuous functions on $\M$.

\begin{definition}
    $\AA$ is the \cs-algebra of continuous sections of $\A_{\hat\I}$.  
    $\Po : \AA \onto \C(\M)$ is the evaluation homomorphism at 
    $\infty\in\hat\I$.
\end{definition}

In most of this paper the quantization map is a map $Q:\C(\M)\to\AA$; 
more generally, the domain of $Q$ may be only a dense subalgebra of 
$\C(\M)$, but it must contain the smooth functions $\C^\infty(\M)$.  
The composition $\Po\circ Q$ is required to be the identity map; that 
is, applying the quantization map to a function $f\in\C(\M)$ gives a 
continuous section of $\A_{\hat\I}$ whose value at $\infty\in\hat\I$ 
is $f$.  Finally, $Q$ is required to commute with the involution 
($*$-structure) and fit a relation with the Poisson bracket.

Specifying the total quantization map, $Q$, is equivalent to 
specifying, for each point $i\in\I:=\hat\I\smallsetminus\{\infty\}$, a 
map $Q_i$ to $\A_i$ which is just $Q$ composed with the evaluation at 
$i$.  It may be possible to reconstruct the continuous field 
$\A_{\hat\I}$ from this system of quantization maps.  The total 
quantization map $Q$ can be reconstructed as the direct product of the 
$Q_i$'s.  The codomain of this reconstructed quantization map is 
superficially the \cs-algebraic direct product $\prod_{i\in\I} \A_i$; 
however, the image of $Q$ is actually contained in 
$\AA\subset\prod_{i\in\I} \A_i$ and will usually generate $\AA$ (as a 
\cs-algebra) in most cases.  Given a proposed system of quantization 
maps $Q_i$, it is a nontrivial convergence condition that the image of 
$Q$ consists of sections of some continuous field.

For some purposes, including defining quantization of a vector bundle, 
the continuous field $\A_{\hat\I}$ is the only important structure.  
For this reason, in \cite{haw1}, I gave the following:
\begin{definition}
A \emph{general quantization} of $\M$ is a continuous field 
$\A_{\hat\I}$ with $\C(\M)$ as the fiber over $\infty$.
\end{definition}
This is just enough structure to define whether or not a 
sequence of operators converges to a given function on $\M$.

The context of this paper is the geometric quantization of compact 
\Kahler\ manifolds.  In this case $\hat\I=\hat\N := 
\{1,2,\ldots,\infty\}$, the 1-point compactification of the positive 
integers.  The algebras $\A_N$ for each $N\in\N$ are 
(finite-dimensional) matrix algebras, and if $\M$ is connected, they are 
simple (i.~e., ``full'') matrix algebras.

Assume for simplicity that $\M$ is connected.  then, all the information of 
the general quantization is contained in $\AA$ and $\Po$.  The index 
set $\hat\N$ can be recovered as the spectrum of the center of $\AA$.

The ideal $\AA_0 := \ker \Po$ is the algebra of sections of 
$\A_{\hat\N}$ that vanish at $\infty$.  Equivalently, $\AA_0$ consists 
of those sections over $\N$ for which the sequence of norms converges 
to $0$.  Since $\N$ is discrete, $\AA_0$ is just the \cs-algebraic 
direct sum $\bigoplus_{N\in\N} \A_N$.

I shall be concerned with two choices of quantization maps here.  The 
first are the Toeplitz quantization maps $T_N$; these are manifestly 
completely positive and thus defined on all of $\C(\M)$.  The second 
type are the geometric quantization maps $\gq_N$; as I shall show in 
Thm.~\ref{GQ.equivalent}, these correspond to the same general 
quantization as the Toeplitz quantization maps do.

Note that a general quantization can be phrased as an extension 
\beq 
0\to\AA_0\longrightarrow\AA\stackrel{\Po}\longrightarrow\C(\M)\to0 
\mbox.\label{ext}
\eeq 
The total Toeplitz quantization map $T : \C(\M)\to\AA$ gives a 
completely positive splitting of \eqref{ext}.

\subsection{Quantized vector bundles}
The category equivalence of vector bundles over $\M$ with finitely 
generated projective (f.~g.~p.)  modules of $\C(\M)$ is well known.  
The $\C(\M)$-module corresponding to a vector bundle $V$ over $\M$ is 
the space of continuous sections $\Gamma(\M,V)$.  This suggests the 
following definition (see \cite{haw1}).
\label{V.def}\begin{definition} 
Given a general quantization, expressed as $\Po:\AA\onto\C(\M)$, a 
quantization of a vector bundle $V$ is any f.~g.~p.\ $\AA$-module, 
$\VV$, such that the push-forward by $\Po$ is 
$\Po_*(\VV)=\Gamma(\M,V)$.
\end{definition}

For every $i\in\hat\I$, pushing $\VV$ forward by the evaluation homomorphism 
gives a module $V_i$ of $\A_i$.  The $\AA$-module $\VV$ is equivalent 
to a bundle of modules over $\hat\I$ whose fiber over $i$ is $V_i$.

It is not obvious \emph{a priori} that any quantization of $V$ will 
exist, or that it will be at all unique.  To investigate these issues, 
it is helpful to consider $K$-theory; the group $K^0(\M)$ classifies 
vector bundles; the group $K_0(\AA)$ classifies f.~g.~p.\ $\AA$-modules 
(which are the quantized vector bundles).

The short exact sequence \eqref{ext} leads, as usual, to a six-term, 
periodic exact sequence in $K$-theory; incorporating the identity 
$K_*[\C(\M)] = K^*(\M)$, this reads,
\beq
\begin{CD}
K_0(\AA_0) @>{\beta}>> K_0(\AA) @>>> K^0(\M) \\
@A{\alpha}AA @. @VVV \\
K^1(\M) @<<< K_1(\AA) @<<< K_1(\AA_0)\mbox.\\
\end{CD}
\label{periodic.exact}\eeq

Assume that $\I$ is discrete, and the $\A_i$'s are full matrix 
algebras.  Then $\AA_0$ is just the direct sum $\bigoplus_i \A_i$, and 
the $K$-groups are direct sums $K_*(\AA_0) = \bigoplus_i K_*(\A_i)$.  
Since each $\A_i$ is a full matrix algebra, its $K$-theory is 
$K_0(\A_i)\cong\Z$ and $K_1(\A_i)=0$.  This gives that 
$K_0(\AA_0)\cong \Z^{\oplus\infty}$ and $K_1(\AA_0) = 0$.

The algebraic direct sum $\Z^{\oplus\infty}$ is the set of sequences 
of integers with finitely many nonzero terms.  This corresponds to the 
fact that an f.~g.~p.\ module of $\AA_0$ is a direct sum of (finitely 
generated) modules of finitely many $\A_N$'s.  Any $\AA_0$-module is 
also an $\AA$-module; the f.~g.~p.\ $\AA_0$-modules are precisely 
those $\AA$-modules which are finite-dimensional as vector spaces.  
This identification corresponds to the map $\beta$ in 
\eqref{periodic.exact}, and shows that $\beta$ must be injective, 
which, by exactness, shows that the map $\alpha$ is zero.  With this, 
the exact sequence \eqref{periodic.exact} breaks down into the 
isomorphism
\[
K_1(\Po): K_1(\AA) \isom K^1(\M)
\]
and the short exact sequence
\beq
0 \to \Z^{\oplus\infty} \longrightarrow K_0(\AA) 
\xrightarrow{K_0(\Po)} K^0(\M) \to 0
\mbox.\label{short.exact}\eeq

This $K_0(\Po)$ maps the $K$-class of a quantization of a vector 
bundle, $V$, to the $K$-class of $V$.  The exact sequence 
\eqref{short.exact} shows that $K_0(\Po)$ is surjective, which 
suggests that any vector bundle can be quantized.
 \begin{thm}
\label{existence}
For a general quantization $\A_{\hat\I}$ such that $\I$ is discrete 
and the fibers over $\I$ are matrix algebras, every vector bundle can 
be quantized.
\end{thm}
 \begin{proof}
Any vector bundle can be realized as the image of some idempotent 
matrix of continuous functions $e\in M_m[\C(\M)]$.  The ideal $\AA_0$ 
consists of all compact operators in $\AA$, so any element of the 
preimage $\Po^{-1}(e)$ is ``essentially'' idempotent (i.~e., modulo 
compacts).  There therefore exists an idempotent $\tilde e$ such that 
$\Po(\tilde e) = e$.  The right image $\AA^m \tilde e$ is a 
quantization of the vector bundle in question.
\end{proof}

The short exact sequence \eqref{short.exact} also shows that $\ker 
K_0(\Po) = \Z^{\oplus\infty}$. This means that the $K$-class of a 
quantization of $V$ is uniquely determined by $V\!$, modulo 
$\Z^{\oplus\infty}$. This suggests:
 \begin{thm}
\label{uniqueness}
With the hypothesis of Theorem~\ref{existence}, quantization of 
a vector bundle is unique modulo finite-dimensional modules.
\end{thm}
 \begin{proof}
We need to prove that if $\VV$ and $\VV'$ are quantizations of 
$V\!$, then there exists a module homomorphism from $\VV$ to 
$\VV'$ whose kernel and cokernel are finite-dimensional (in other 
words, a Fredholm homomorphism).

Any f.~g.~p.\ module can be realized as the (right) image of an 
idempotent matrix over $\AA$.  So, identify $\VV$ and $\VV'$ with the 
images of idempotents $e,e'\in M_m(\AA)$; that is, $\VV=\AA^m e$ and 
$\VV'=\AA^m e'\!$.  These idempotents can be chosen so that 
$\Po(e)=\Po(e')$; therefore, $e-e'\in M_m(\AA_0)$.  Multiplication by 
$e'$ (respectively, $e$) gives a homomorphism $\varphi' : \VV 
\xrightarrow{\cdot e^{\smash\prime}} \VV'$ (resp., $\varphi :\VV' 
\xrightarrow{\cdot e} \VV$).

Let $k$ be the self-adjoint idempotent whose (right) image is $\ker 
\varphi\circ\varphi'\!$.  Since this is a subspace of $\VV\!$, $k$ 
satisfies $ke=k$, and since this is the kernel of 
$\varphi\circ\varphi'\!$, $k$ satisfies $ke'e=0$.  \emph{A priori} $k$ 
is not necessarily in $M_m(\AA)$.  However, the entries of $k$ are in the 
\cs-algebraic direct product of the $\A_i$'s; i.~e., bounded sections 
of $\A_{\hat\I}$ over $\I$.  $\AA_0$ is an ideal in this algebra, so 
$k = k (e-e')e \in M_m(\AA_0)$.  Therefore, $\ker 
\varphi\circ\varphi'$ is an f.~g.~p.\ module of $\AA_0$ and is thus 
finite-dimensional.  By an identical argument, $\ker 
\varphi'\circ\varphi$ is also finite-dimensional.  This implies that 
the kernel and cokernel of $\varphi'$ are finite-dimensional.
\end{proof}

The converse is clearly also true: If $\VV$ is a quantization of 
$V\!$, and $\VV'$ is isomorphic to $\VV$ modulo finite-dimensional 
modules, then $\VV'$ is also a quantization of $V\!$.

\section{Toeplitz Quantization}
\label{Toplitz}
Again, let $\M$ be a compact, connected \Kahler\ manifold. Now, let $L$ be 
a Hermitian line bundle with curvature given by the symplectic form as
$\nabla^2=-i\omega$ and $L_0$ a holomorphic line bundle with an inner 
product on sections (i.~e., a pre-Hilbert structure on 
$\Gamma(\M,L_0)$).
\begin{definition}
For each $N\in\N$, $\LN:=L_0\otimes L^{\otimes N}\!$, $\HN := 
\Gh(\M,\LN)$ 
(the space of holomorphic sections of $\LN$), and $\A_N := \End \HN$ 
(matrices over $\HN$). 
\end{definition}

The inner products on sections of $L_0$ and fibers of $L$ combine to 
give an inner product on sections of $\LN$.  This makes $\HN$ into a 
Hilbert space; it does not need to be completed, since it is 
finite-dimensional.  The Hilbert space structure of $\HN$ makes $\A_N$ 
a \cs-algebra.

The connections and inner products must be compatible.  For 
convenience, assume that $L_1$ (and thus any $L_{N\in\N}$) is ``positive'' 
(if not, just reparameterize $N$).  This guarantees that $\A_N$ is 
nontrivial for all $N\in\N$.  The simplest choice of $L_0$ is just the 
trivial line bundle with the trivial connection and the inner product 
given by integrating with the canonical volume form $\omega^n/n!$ 
($n:=\dim_\co \M$).

The space $\HN$ is naturally a Hilbert subspace of $L^2(\M,\LN)$ which 
is a subspace of the the Hilbert space of $(0,*)$-forms with 
coefficients in $\LN$. 
\begin{definition}
Let $\Pi_N$ be the self-adjoint projection onto $\HN$.  The Toep\-litz 
quantization map $T_N: \C(\M) \to \A_N$ is given by \beq T_N(f) := 
\Pi_N f \mbox.%\label{PNdef}
\eeq
\end{definition}
In other words, the action of $T_N(f)$ on an element of $\HN$ 
(holomorphic section of $L_N$) is given by first multiplying by $f$ 
(giving a non-holomorphic section) and then projecting back to $\HN$ 
by $\Pi_N$.  This $T_N$ is automatically a unital and (completely) positive 
map; therefore, it is norm-contracting.

\subsection{Vector bundles}
Suppose that we are given a smooth vector bundle $V$ with a specific 
connection.  I would like to construct from $V$ a sequence of $\A_N$ 
modules.  The algebra $\A_N$ can be written as $\A_N = \End \HN = 
\Hom(\HN,\HN)$ and can be thought of as consisting of square matrices 
of height and width $\HN$.  Any module of $\A_N$ can be written as 
$\Hom(E,\HN)$ and thought of as consisting of rectangular matrices of 
height $\HN$ and width $E$.  Any construction for $E$ should 
generalize that of $\HN$.

Thanks to the Kodaira vanishing theorem (see the appendix and e.~g., 
\cite{g-h1}) and the assumption that $L_1$ is positive, 
$\HN=\Gh(\M,\LN)$ can also be realized as the kernel of the 
$\LN$-twisted Dolbeault operator that acts on $\Omega^{0,*}(\M,\LN)$.  

In order to generalize $\HN$ appropriately, we will need:
\begin{definition}
$D_V:=\nabla_{\!\bar\partial}+(\nabla_{\!\bar\partial})^*=i\gamma^\mu 
\nabla_{\!\mu}$ is the $V^*\otimes\LN$-twisted Dolbeault operator, a 
Dirac-type operator acting on the smooth $(0,*)$-forms 
$\Omega^{0,*}(\M,V^*\otimes\LN)$.  
\end{definition} 
Here, $\nabla$ is the connection, and the Dirac matrices satisfy 
$[\gamma^\mu,\gamma^\nu]_+ = \gamma^\mu \gamma^\nu + \gamma^\nu 
\gamma^\mu = g^{\mu\nu}\!$, which differs from the usual convention by 
a factor of 2.  
A number of inequalities related to Dolbeault operators will prove 
useful, but the proofs of these are relegated to the appendix.

Natural generalizations of $\A_N$ and $T_N$ are,
\begin{definition}
$\tilde V_N := \Hom(\tilde\lef^V,\HN)$ where $\tilde\lef^V=\ker D_V$.  
The map $T_N^V:\Gamma(\M,V)\to\tilde V_N$ is given by
\beq 
T_N^V(v) := \Pi_N v \mbox.\label{PVdef}
\eeq
\end{definition}
In this, multiplication must be understood to mean contraction of $V$ with 
$V^*\!$.  Multiplying an element of 
$\tilde\lef^V\subset\Omega^{0,*}(\M,V^*\otimes\LN)$ by $v$ gives an 
element of\fchoice{\linebreak}{} $\Omega^{0,*}(\M,\LN)$; $\Pi_N$ then projects this down to 
$\HN$.  If $V$ is trivial (i.~e., $V=\co\times\M$ with the trivial 
connection), then $T_N^V$ reduces to $T_N$.

The tildes will be dispensed with by the end of the next section.

\section{Convergence}
\label{convergence}
If $V$ has an inner product --- or if we assign one --- then there is 
a natural operator norm on $\tilde{V}_N$ which generalizes the norm on 
$\A_N$.  This corresponds to the norm $\norm{v} := \sup_{x\in\M} 
\norm{v(x)}$ on sections.  None of the \emph{constructions} here will 
require an inner product on $V$; however, several of the \emph{proofs} 
will make use of one --- which can be taken arbitrarily.

The main property of the maps $T_N$ and $T_N^V$ that is needed to 
prove convergence of the quantization of both the algebra and vector 
bundles is
 \begin{lem}\label{converge.lemma}
For any function $f\in\C(\M)$ and any section $v\in\Gamma(\M,V)$, 
\[
\lim_{N\to\infty} \left\Vert T_N(f)T_N^V(v) - 
T_N^V(fv)\right\Vert = 0
\mbox.\]
\end{lem}
 \begin{proof}
Let $D$ be the $\LN$-twisted Dolbeault operator, so that $\HN=\ker D$.  
We can approximate $\Pi_N$ by $(1 + \alpha D^2)^{-1}$ with $\alpha$ a 
positive real number; in fact, by Lemma~\ref{spectrum.lemma} (in the 
appendix) there exists a constant $C<1$ such that \beq \spec D^2 
\subset \{0\}\cup[N-C,\infty) \mbox,\label{spectrum2}\eeq so the error 
in $\Pi_N\approx(1+\alpha D^2)^{-1}$ is bounded as
\[
\Norm{\Pi_N - (1 + \alpha D^2)^{-1}} \leq (1 + \alpha [N-C])^{-1} 
\leq \alpha^{-1} (N-C)^{-1}
\mbox.\] 

Now, for $f\in\C^\infty(\M)$ any smooth function, the commutator 
$[f,\Pi_N]_-$ is approximated by
\begin{align*}
[f,(1+\alpha D^2)^{-1}]_- &= \alpha (1+\alpha D^2)^{-1} [D^2,f]_- 
(1+\alpha D^2)^{-1} 
\\&= i\alpha (1+\alpha D^2)^{-1} [D,\gamma^\mu f_{|\mu}]_+ (1+\alpha 
D^2)^{-1}
\mbox.\end{align*}
Another consideration of \eqref{spectrum2} shows that,
\[
\left\Vert(1+\alpha D^2)^{-1} D\right\Vert \leq 
\frac{(N-C)^{1/2}}{1+\alpha (N-C)} 
\leq \alpha^{-1} (N-C)^{-1/2} 
\mbox.\] 
So,
\[
\left\Vert[f,(1+\alpha D^2)^{-1}]_-\right\Vert \leq 2 (N-C)^{-1/2} 
\left\Vert\gamma^\mu f_{|\mu}\right\Vert = \sqrt{2}\, (N-C)^{-1/2} 
\norm{\nabla f}
\mbox.\]
This gives that
\[
\left\Vert[f,\Pi_N]_-\right\Vert \leq \sqrt{2}\, (N-C)^{-1/2} 
\norm{\nabla f} + 2 \alpha^{-1} (N-C)^{-1} \norm{f}
\mbox,\]
for \emph{any} $\alpha>0$, and therefore,
\[
\left\Vert[f,\Pi_N]_-\right\Vert \leq \sqrt{2} \, (N-C)^{-1/2} 
\norm{\nabla f}
\mbox.\]

With a slight abuse of notation,
\[
T_N(f) T_N^V(v) - T_N(fv) = T_N^V(f\,\Pi_N v - fv) = 
T_N^V([f,\Pi_N]_-v)
\mbox.\]
By construction, $T_N^V$ is norm-contracting; thus, 
\[
\left\Vert T_N(f) T_N^V(v) - T_N^V(fv)\right\Vert \leq 
\sqrt{2}\, (N-C)^{-1/2} \norm{\nabla f}\, \norm{v} 
\mbox.
\] 
So $\left\Vert T_N(f) T_N^V(v) - T_N^V(fv)\right\Vert \to 0$ as 
$N\to\infty$.  Since $T_N$ and $T_N^V$ are norm-contracting, they are 
continuous, and since $\C^\infty(\M)\subset\C(\M)$ is a dense 
subalgebra, the conclusion holds for all $f\in\C(\M)$.
\end{proof}

\begin{definition}
The total Toeplitz quantization map $T:\C(\M)\to \prod_{N\in\N} \A_N$ 
is the direct product of the $T_N$'s.  Also 
$\AA_0:=\bigoplus_{N\in\N}\A_N$ is the \cs-algebraic direct sum, and 
$\AA:=\Im T + \AA_0$.
\end{definition}
 \begin{lem}
$\AA$ is a \cs-algebra, and $T$ induces an isomorphism 
$\C(\M)\isom\AA/\AA_0$.
\end{lem}
 \begin{proof}
Lemma~\ref{converge.lemma} is in this case equivalent to the 
statement that for any functions $f,g\in\C(\M)$,
\beq
T(f) T(g) - T(fg) \in \AA_0
\mbox.\label{modA0}\eeq
The direct sum $\AA_0$ is an ideal in the direct product 
$\prod_{N\in\N} \A_N$, so \eqref{modA0} shows that $\AA$ is 
algebraically closed.  Since $T$ is norm-contracting, $\Im T$ is 
norm-closed, and so $\AA$ is norm-closed.  Hence, $\AA$ is a 
\cs-algebra.

Equation \eqref{modA0} also shows that $T$ induces (by composition 
with the quotient map $\AA\onto\AA/\AA_0$) a $*$-homomorphism 
$\C(\M)\onto\AA/\AA_0$.  This is surjective because of the definition 
of $\AA$.  We need to verify that it is injective.

Since $\AA$ lies inside the direct product of the $\A_N$'s, there is 
for each $N\in\N$ an obvious ``evaluation'' homomorphism 
$\Po_N:\AA\onto\A_N$.  Define the normalized partial traces $\ntr_N : 
\AA \to \co$ by $\ntr_N a := \tr[\Po_N (a)]/\!\dim \HN$, so that 
$\ntr_N 1 = 1$.  The normalized trace is norm-contracting, so any 
$a\in\AA$ satisfies $\lvert \ntr_N a\rvert \leq \norm{\Po_N(a)}$; 
therefore,
\beq
a\in\AA_0 \implies 
\lim_{N\to\infty} \ntr_N a = 0 
\mbox.\label{ntrA0}\eeq

Note that for any $f\in\C(\M)$, the (unnormalized) trace of $T_N(f)$ 
can be expressed as
\[
\tr[T_N(f)] = \Tr [\Pi_N f] = \lim_{t\to\infty} \Tr[e^{-tD^2}f] 
\mbox.\] 
Using the asymptotic expansion for $e^{-tD^2}$ (the ``heat kernel 
expansion''\!, see \cite{gil1}), this can be evaluated explicitly as a 
polynomial in the curvatures of $\TM$ and $\LN$.  The result is a 
polynomial in $N$ with leading order term
\[
\left(\tfrac{N}{2\pi}\right)^n\intM f\frac{\omega^n}{n!}
\mbox.\]

This, with \eqref{ntrA0}, shows that $\ntr_\infty := \lim_{N\to\infty} 
\ntr_N$ is well-defined on $\AA$, vanishes on $\AA_0$ and satisfies
\[
\ntr_\infty [T(f)] = \frac1{\vol \M} \intM f \frac{\omega^n}{n!} 
\mbox.\] 
Suppose that some function $f$ is in the kernel of the induced 
homomorphism \fchoice{\linebreak}{} $\C(\M)\to\AA/\AA_0$, or 
equivalently that $T(f)\in\AA_0$.  The kernel of a $*$-homomorphism is 
spanned by its positive elements, so we can assume without loss of 
generality that $f\geq 0$.  This implies that $0= \ntr_\infty[T(f)] 
\propto \intM f \omega^n/n!$, but since $\omega^n$ is nonvanishing 
this implies $f=0$.  So, the homomorphism is injective and thus an 
isomorphism.
\end{proof}

\begin{definition}
$\Po:\AA\onto\C(\M)$ is the composition of the natural surjection 
$\AA\onto\AA/\AA_0$ with the inverse of the isomorphism induced by $T$.
\end{definition}

The following shows that $\AA$ indeed gives a general quantization of 
$\M$.
 \begin{thm}\label{converge}
There is a continuous field of \cs-algebras over $\hat\N$ such that 
the fiber over $N\in\N$ is $\A_N$, the fiber over $\infty$ is 
$\C(\M)$, and the algebra of continuous sections is $\AA$.
\end{thm}
 \begin{proof}
Let $\Po_N:\AA\to\A_N$ denote the evaluation map at $N$.  Most of the 
axioms given in \cite{dix1} for a continuous field of \cs-algebras are 
easily verified.  The nontrivial axiom is the requirement that for any 
$a\in\AA$ the norms $\norm{\Po_N(a)}$ define a continuous function on 
$\hat\N$.  Since continuity is only an issue at $\infty\in\hat\N$, 
this reduces to the requirement that $\norm{\Po_N(a)}\to\norm{\Po(a)}$ 
when $N\to\infty$.  It is sufficient to prove this on $\Im T$; in 
other words, we need to show that for any $f\in\C(\M)$,
\[
\lim_{N\to\infty} \norm{T_N(f)} = \norm{f}
\mbox.\]

The spectrum $\hat\AA$ (of irreducible representations, see 3.2.2 of 
\cite{dix1}) is a non-Haus\-dorff union of $\M$ and $\N$, although it 
maps continuously onto $\hat\N$.  According to Prop.~3.3.2 of 
\cite{dix1}, the function on $\hat\AA$ defined by the norms of the 
images of any $a\in\AA$ is lower semi-continuous.  This means that for 
any $x\in\M$,
\begin{align*}
\liminf_{N\to\infty} \norm{T_N(f)} &\geq \norm{f(x)} \\
&\geq\norm{f}
\mbox.\end{align*}
On the other hand, because $T_N$ is norm-contracting,
\[
\limsup_{N\to\infty} \norm{T_N(f)} \leq \sup_{N\in\N} 
\norm{T_N(f)} 
\leq \norm{f}
\mbox.\]
\end{proof}

Each of the maps $T_N$ is surjective (Prop.~4.2 of \cite{b-m-s}), so 
this is clearly the smallest continuous field such that $N\mapsto 
T_N(f)$ defines a continuous section.  In fact, $\Im T$ generates 
$\AA$ as a \cs-algebra.

\begin{definition}\label{VV.def}
\[
T^V:\Gamma(\M,V)\to \prod_{N\in\N} \tilde V_N
\]
is the direct product of the $T^V_N$'s, $\VV:=\AA \cdot \Im T^V$ is 
the $\AA$-module generated by $\Im T^V\!$, $V_N$ is the restriction of 
$\VV$ to an $\A_N$-module, and $\Po^\VV : 
\VV\onto\VV/\AA_0\!\VV=\Po_*(\VV)$ is the natural surjection.
\end{definition}
The following lemma shows that the analytic condition $\norm{v_N}\to0$ 
can be expressed algebraically.
 \begin{lem}\label{VV.preliminary}
\[
\AA_0\!\VV = \Bigl\{v\in\prod_{N\in\N}V_N \Big| 
\lim_{N\to\infty} \norm{v_N} = 0 \Bigr\}
\mbox,\] 
and $T^V$ induces a homomorphism of 
$\C(\M)$-modules, $\Po^\VV\circ T^V:\Gamma(\M,V)\onto\Po_*(\VV)$.
\end{lem}
 \begin{proof}
As I have mentioned, $T_N^V$ is norm-contracting; thus $T^V(v)$ is 
bounded.  Because of this, the sequence of norms coming from any 
element of $\AA_0\!\VV=\AA_0\Im T^V$ must converge to $0$.  
Conversely, $\AA_0\!\VV$ is norm-closed and contains all sequences in 
$\prod_{N\in\N}V_N$ with finitely many nonzero terms.  This proves the 
first claim.

With this, Lemma~\ref{converge.lemma} shows 
that for all $f\in\C(\M)$ and $v\in\Gamma(\M,V)$
\[
T(f) T^V(v) - T^V(fv) \in \AA_0\!\VV
\mbox,\]
which proves the second claim.
\end{proof}

 \begin{lem}\label{generation}
For $N$ sufficiently large, $V_N=\tilde V_N$. 
\end{lem}
 \begin{proof}
Equivalently, for sufficiently large $N$, the image of $T_N^V$ 
generates $\tilde{V}_N$ as an $\A_N$-module.  If not, then $\Im 
T_N^V$ must lie inside a proper submodule of $\tilde{V}_N$, and so there 
must exist $\psi\in\tilde\lef^V$ such that, for any $\varphi\in\HN$ and 
$v\in\Gamma(\M,V)$, $\langle\varphi\rvert T_N^V(v) \lvert\psi\rangle 
= 0$.

Take any nonzero $\varphi\in\HN$ and $\psi\in\tilde\lef^V$.  Let 
$\psi_0\in\Gamma^\infty(\M,V^*\otimes\LN)$ be the component of $\psi$ 
in degree $0$.  Assume $N$ to be sufficiently large that $\psi_0$ is 
guaranteed by Corollary~\ref{vanish2} not to vanish.  Using the fact 
that $\M$ is connected, the zeros of $\varphi$ must form a proper 
subvariety of $\M$, and $\psi_0$ must be nonzero on an open set; 
therefore, there exists $y\in\M$ where $\varphi(y),\psi_0(y)\neq0$.  
If $v\in\Gamma(\M,V)$ approximates the distribution 
$\varphi(y)\bar\psi_0(y)\delta(x,y)$, then $\langle\varphi\rvert 
T_N(v) \lvert\psi\rangle = \langle\varphi\rvert v \lvert\psi\rangle$ 
will approximate $\norm{\varphi(y)}^2 \norm{\psi_0(y)}^2$ and must be 
nonzero for a sufficiently close approximation.
\end{proof}

\subsection{Category}
It remains to be proven that $\VV$ is a finitely generated, projective 
(f.~g.~p.)  module and that its push-forward is 
$\Po_*(\VV)=\Gamma(\M,V)$.  To do this, it will be helpful to make the 
correspondence $V\mapsto\VV$ into a functor.  Since the module $\VV$ 
is not constructed from the vector bundle $V$ alone but from $V$ 
accompanied by a connection, the domain of this functor must be a 
category of vector bundles with connections.  We need to identify 
those bundle homomorphisms which will lead naturally to module 
homomorphisms.

A bundle homomorphism naturally defines a map of sections.  It also 
naturally gives (by tensor product with the identity map) 
homomorphisms of the tensor products with any other bundle (such as 
1-forms).  For simplicity, I will denote all these trivially derived 
maps by the same symbol as the original homomorphism.
\begin{definition}
a morphism of bundles with connections, $\phi: V\to W$, is a smooth 
bundle homomorphism such that for any smooth section 
$v\in\Gamma^\infty(\M,V)$, 
\beq
\label{morph.def} 
\phi(\nabla_{\!  V} v) = \nabla_{\!  W} \phi(v) 
\mbox;\eeq 
in other words, $\phi$ is covariantly constant.
\end{definition}
With these morphisms, vector bundles with connections form an Abelian 
category.  Clearly, the identity homomorphism on any bundle satisfies 
the above property, and the composition of two such morphisms does as 
well.  Also, the kernel and cokernel of such a morphism inherit 
natural connections.

Now, let's try and construct a functor $\func$ from this category of 
bundles with connections to the category of $\AA$-modules, such that 
$\func(V)=\VV\!$.

\begin{dthm}
    Any morphism $\phi:V\to W$ of bundles with connections induces a 
    homomorphism $\func(\phi):\VV\to\WW$ of $\AA$-modules which 
    satisfies,
    \beq
    T^W\circ\phi = \func(\phi)\circ T^V
    \mbox.\label{T.func}\eeq
\end{dthm}
\begin{proof}
$\phi$ gives an adjoint map on the dual bundles in the opposite 
direction, and in turn maps the spaces of forms 
$\phi^*:\Omega^{0,*}(\M,W^*\otimes\LN)\to\Omega^{0,*}(\M,V^*\otimes\LN)$. 
Because $\phi$ intertwines connections, the map $\phi^*$ intertwines 
Dolbeault operators.  If $\psi\in\ker D_W = \tilde\lef^W$, then 
$D_V\phi^*(\psi)=\phi^*(D_W\psi)=0$; this means that the restriction 
of $\phi^*$ to $\tilde\lef^W$ maps 
$\phi^*:\tilde\lef^W\to\tilde\lef^V$.  This induces a homomorphism 
$\phi_*:\tilde V_N\supset V_N\to\tilde W_N$.

Put these maps together to define $\func(\phi)$.  \emph{A priori}, 
$\func(\phi)$ maps an element of $\VV$ to some sequence of elements of 
$\tilde W_N$.  We need to prove that the image of $\func(\phi)$ in 
fact lies inside $\WW$.

For any $v\in\Gamma(\M,V)$ and $\psi\in\tilde\lef^W$,
\[
T_N^W[\phi(v)]\psi=\Pi_N\phi(v)\psi = \Pi_N v\phi^*(\psi) = 
\phi_*\!\left[T_N^V(v)\right]\psi
\mbox.\]
So, $T_N^W\circ\phi = \phi_*\circ T_N^V$ and hence, 
Eq.~\eqref{T.func}.
Since 
$\WW$ is defined to be generated by $\Im T^W\!$, this shows that 
indeed $\func(\phi):\VV\to\WW$.
\end{proof}

$\func$ is an additive functor.  It respects identity maps, 
compositions, and sums of morphisms.  Because of this, $\func$ must 
respect finite direct sums; this property is also easily seen from the 
construction of $\func$.

The category of bundles with connections actually behaves somewhat 
trivially.  Because a morphism is covariantly constant, it can be 
specified completely by its action at a single point.  As a result, 
this category behaves somewhat like the category of finite-dimensional 
vector spaces.  Any short exact sequence splits.  Because of this, any 
additive functor (such as $\func$) on this category is exact.

Of course, not all bundle homomorphisms are morphisms of bundles with 
connections.  We will need some module homomorphisms that do not come 
from $\func$.  The following result shows that an isomorphism of 
vector bundles can be used to construct a homomorphism of modules 
which is an isomorphism modulo finite-dimensional modules.

 \begin{lem}\label{almost.isomorphic}
Let $V$ and $W$ be isomorphic bundles with different connections.  
Then there exists a Fredholm homomorphism $u:\VV\to\WW$ (compare 
Thm.~\ref{uniqueness}), which 
satisfies 
\beq
\Po^\WW\circ u\circ T^V =\Po^\WW\circ T^W
\label{T.u}\mbox.\eeq
\end{lem}
\begin{proof}
The homomorphism $u$ is specified by giving, for each $N$, a 
homomorphism $u_N:V_N\to W_N$ of $\A_N$-modules.  Define $\Pi_N^V$ to 
be the spectral projection at $0$ for the $V^*\otimes\LN$-twisted 
Dolbeault operator (likewise with $W$); that is, $\Pi_N^V$ is an 
idempotent with $\Im\Pi_N^V=\ker D_V$ and $\ker\Pi_N^V=\Im D_V$.

The isomorphism of $V$ and $W$ gives a natural (isometric) inclusion 
$\iota:\tilde\lef^W\into\Omega^{0,*}(\M,V^*\otimes\LN)$.  Composing 
this with $\Pi_N^V$ gives $\Pi_N^V\iota:\tilde \lef^W\to\tilde\lef^V$, 
and $\Pi_N^W\iota$ is the identity on $\tilde\lef^W$.  According to 
Lemma~\ref{Pi.norms},
$$
\lim_{N\to\infty} \norm{\Pi_N^V-\Pi_N^W} = 0
\mbox;\eqno{\eqref{norm.Pi.diff}}
$$
therefore, for $N$ sufficiently large, $\norm{\Pi_N^V-\Pi_N^W}<1$. 
When this is so, $\Pi_N^V\iota$ is injective, because if 
$\psi\in\tilde\lef^W$ is nonzero then,
\begin{align*}
\Norm{\Pi_N^V\iota\psi} &= \Norm{(\Pi_N^V-\Pi_N^W)\iota\psi+\psi}\\
&\geq \left(1-\Norm{\Pi_N^V-\Pi_N^W}\right)\norm\psi > 0
\mbox.\end{align*}
The existence of a similar injection in the opposite direction 
establishes that $\Pi_N^V\iota$ is bijective.

Recall from Lemma~\ref{generation} that for $N$ sufficiently large 
$V_N=\tilde{V}_N$ (and likewise with $W$).  When $N$ is sufficiently 
large, we can define $u_N:V_N\to W_N$ to be the bijection given by 
$u_N(v_N)=v_N\Pi_N^V\iota$.  For small $N$, it doesn't matter what 
$u_N$ is.

Now assemble the $u_N$'s into $u$. The kernel and cokernel of $u$ come 
entirely from the finitely many $u_N$'s which are not bijective, and 
thus are finite-dimensional. In other words, $u$ is Fredholm.

But does the image in fact lie inside $\WW$?  Using 
Eq.~\eqref{norm.Pi.diff} again shows that, for any $v\in\Gamma(\M,V)$,
\[
\Norm{T_N^V(v)\Pi_N^V\iota - T_N^W(v)} 
\leq \Norm{T_N^V(v)}\Norm{\Pi_N^V-\Pi_N^W} \to0
\]
as $N\to\infty$; therefore, by Lemma~\ref{VV.preliminary}, 
$u[T^V(v)]-T^W(v)\in \AA_0\!\WW\subset\WW$.
\end{proof}

 \begin{thm}\label{fgp.theorem}
$\VV=\func(V)$ is a 
quantization of $V$ by the definition in Sec.~\ref{V.def}.
\end{thm}
 \begin{proof}
For any vector bundle $V\!$, there exists another vector bundle $W$ 
such that the direct sum is some trivial bundle $V\oplus W \cong 
\co^m\times\M$. Choose an arbitrary connection on $W$. 
As noted above, $\func$ respects finite direct sums, so
$\func(V\oplus W)=\VV\oplus\WW$.

By Lemma~\ref{almost.isomorphic} there exists an $\AA$-module 
homomorphism $u : \VV\oplus\WW \to \AA^m$ whose kernel and cokernel 
are finite-dimensional and thus projective.  All the terms of the 
exact sequence
\[
0 \to \ker u \longrightarrow \VV\oplus\WW 
\stackrel{u}{\longrightarrow} \AA^m \longrightarrow \coker u \to 0
\mbox,\]
other than $\VV\oplus\WW$, are now seen to be f.~g.~p.\ modules; 
therefore $\VV\oplus\WW$, and thus $\VV$, is f.~g.~p.

It remains to prove that $\Po_*(\VV)=\VV/\AA_0\!\VV=\Gamma(\M,V)$.  
Lemma~\ref{VV.preliminary} showed that $\Po^\VV\circ 
T^V:\Gamma(\M,V)\onto\Po_*(\VV)$ is a $\C(\M)$-module homomorphism, 
and it is clearly surjective by the definition of $\VV\!$.  We need to 
prove that the kernel of $\Po^\VV\circ T^V$ is trivial.

Let $\phi$ denote the natural inclusion $\phi : V \into V\oplus W$ (as 
bundles with connections) and $\varphi$ the equivalent inclusion 
$\varphi:V\into\co^m\times\M$ (as a bundle).  If 
$v\in\ker[\Po^\VV\circ T^V]$, then $T^V(v)\in\AA_0\!\VV\!$, so 
$\Po\circ u \circ \func(\phi) \circ T^V(v) = 0$, since 
$u\circ\func(\phi)$ is an $\AA$-module homomorphism.  However, by 
Eq's~\eqref{T.func} and \eqref{T.u},
\[
\Po\circ u \circ \func(\phi) 
\circ T^V = 
\Po \circ u \circ T^{V\oplus W} \circ \phi = \Po \circ 
T \circ \varphi = \varphi
\]
which is injective. 
Therefore, $\ker[\Po^\VV\circ T^V]\subseteq\ker \varphi = 0$, and
\[
\Po^\VV\circ T^V : \Gamma(\M,V) \isom \Po_*(\VV)
\]
is indeed an isomorphism.
\end{proof}

\section{The holomorphic case}
Recall that a holomorphic vector bundle is a bundle with a connection 
whose curvature is of type $(1,1)$.
 \begin{thm}
If $V$ is a holomorphic vector bundle, then for all $N\in \N$,
$V_N = \Hom(\lef^V,\HN)$ where
\[
\lef^V=\Gh(\M,V^*\otimes\LN)
\mbox,\]
and $\Gh$ means holomorphic sections.
\end{thm}
 \begin{proof}
This is much the same as the proof of Lemma~\ref{generation}. 

In the holomorphic case, $D^2$ respects the $\Z$-grading, so 
$\tilde\lef^V$ (and thence $V_N$) is $\Z$-graded.  Sections of $V\!$, 
and thus $\Im T_N^V$, are entirely of degree $0$, so $V_N \subseteq 
\Hom(\lef^V,\HN)$.

If the statement were false, then there would exist a nonzero 
$\psi\in\lef^V$ such that for any, $\varphi\in\HN$ and 
$v\in\Gamma(\M,V)$, $\langle\varphi\rvert v \lvert\psi\rangle =0$.  
However, if $\varphi\neq0$ then the zero sets of $\varphi$ and $\psi$ 
will be proper subvarieties of $\M$; therefore, there exists $y\in\M$ 
where $\varphi(y)\neq0$ and $\psi(y)\neq0$.  So, if $v(x)$ 
approximates the distribution $\varphi(y)\bar\psi(y)\delta(x,y)$, then 
$\langle\varphi\rvert v \lvert\psi\rangle$ will approximate 
$\norm{\varphi(y)}^2\norm{\psi(y)}^2$ and thus be nonzero for a 
sufficiently close approximation.
\end{proof}

\section{Geometric Quantization}
\begin{definition}
The standard geometric quantization maps (see \cite{woo1}) $\gq_N : 
\C^\infty(\M) \to \A_N$ are defined on smooth functions by (with a 
slight abuse of notation)
\beq
\gq_N(f) := \Pi_N \left[f - \tfrac{i}{N} \pi^{\mu\nu} f_{|\mu} 
\nabla_{\!\nu}\right]  
= T_N \! \left(f - \tfrac{i}{N} \pi^{\mu\nu} f_{|\mu} 
\nabla_{\!\nu}\right)
\mbox.\label{GQdef}\eeq
Here $\pi$ is the Poisson bivector, defined by 
$\pi^{\mu\nu}\omega_{\lambda\nu} = \delta^\mu_\lambda$, and $\nabla$ 
is again the connection.
\end{definition}
Following $T_N^V$, there is an obvious generalization of $\gq_N$ 
for vector 
bundles. 
\begin{definition}
$\gq_N^V : \Gamma^\infty (\M,V) \to \tilde V_N$ is given by
\[
\gq_N^V(v) := T_N^V \! \left(v - \tfrac{i}{N} \pi^{\mu\nu} v_{|\mu} 
\nabla_{\!\nu} \right)
\mbox.\]
\end{definition}

 \begin{lem}\label{GQ-T.to0}
For any smooth section $v\in \Gamma^\infty (\M,V)$,
\[
\lim_{N\to\infty}\norm{T_N^V(v)-\gq_N^V(v)} = 0
\mbox.\]
\end{lem}
 \begin{proof}
Let $w^\mu$ be any tangent vector with components in $V\!$, and use $D$ to 
denote both the $\LN$-twisted and $V^*\otimes\LN$-twisted Dolbeault 
operators. Then,
\begin{align*}
-i [D,\gamma_\mu w^\mu]_+ &= [\gamma^\nu \nabla_{\!\nu} , \gamma_\mu 
w^\mu]_+ = [\gamma^\nu, \gamma_\mu w^\mu]_+ \nabla_{\!\nu} + 
\gamma^\nu [\nabla_{\!\nu} , \gamma_\mu w^\mu]_- \\&= w^\mu 
\nabla_{\!\mu} + \gamma^\nu \gamma_\mu w^\mu_{\; |\nu} 
\mbox.\end{align*} 
Because the argument of $T_N^V$ acts between the kernels of the 
Dolbeault operators, this gives the identity
\[ 
0 = T_N^V\left([D,\gamma_\mu w^\mu]_+\right) 
= i T_N^V(w^\mu \nabla_{\!\mu} + \gamma^\nu \gamma_\mu w^\mu_{\; 
|\nu})
\mbox.\]
Now, setting $w^\mu = -\frac{i}{N} \pi^{\nu\mu} v_{|\nu}$ gives
\beq
\gq_N^V(v) - T_N^V(v) = \tfrac{i}N T_N^V(\gamma^\nu \gamma_\mu 
\pi^{\lambda \mu} v_{|\lambda\nu})
\mbox.\label{TQdif}\eeq
Since $T_N^V$ is norm-contracting, for any smooth $v$, the norm of the 
difference \eqref{TQdif} converges to $0$ as $N\to\infty$.
\end{proof}

Equation \eqref{TQdif} is related to a formula due to Tuynman 
\cite{tuy2}.  Namely, for any smooth function $f\in\C^\infty(\M)$,
\[
\gq_N (f) = T_N \left[f + \tfrac1{2N} \Delta(f)\right]
\]
where $\Delta=-\nabla^2$ is the scalar Laplacian.
Since $\Delta$ is a positive operator, this shows that $Q_N$, like 
$T_N$, is (completely) positive, which means that, after all, 
$Q_N$ can be uniquely, continuously defined on all of $\C(\M)$.

As with Toeplitz quantization, we can assemble the $\gq_N$'s 
into a direct-product map $\gq : \C^\infty(\M) \to \prod_{N\in\N} 
\A_N$.
 \begin{thm}\label{GQ.equivalent}
$\gq :\C(\M) \to \AA$ and $\Po\circ\gq = \id$.
\end{thm}
 \begin{proof}
By Lemma~\ref{GQ-T.to0}, for any smooth function $f$, $T(f) - 
\gq(f)\in\AA_0$. This shows that $\Po[\gq(f)] = f$.

Since $\AA = \Po^{-1}[\C(\M)]$, this shows that $\Im \gq \subset \AA$.
\end{proof}
This shows that the general quantization constructed by geometric 
quantization is exactly the same as that constructed by Toeplitz 
quantization.

Analogous to the construction of $\VV\!$, define $\VV'$ to be the 
$\AA$-module 
generated by the image of $\gq^V\!$.
 \begin{thm}
This $\VV'$ is a quantization of $V\!$.
\end{thm}
 \begin{proof}
It is sufficient to prove that $\VV'$ is isomorphic to $\VV$ modulo 
finite-dimensional modules.

Choose some set of smooth sections of $V$ such that their images by 
$T^V$ generate $\VV\!$, and hence their images by $T_N^V$ generate 
$V_N$.  For sufficiently large $N$, their images by $\gq_N^V$ will be 
close enough to those by $T_N^V$ to generate $\tilde V_N=V_N$.  
Therefore, for $N$ sufficiently large $V_N'= V_N$.
\end{proof}

\section{Further structures}
Because $\VV$ has been produced constructively from $V$ and its 
connection, essentially any additional structure that is consistent 
with the connection on $V$ should lift to $\VV\!$.  (This is equally 
true for $\VV'\!$.)

If there is a group $G$ acting on $\M$, and $V$ is a $G$-equivariant 
vector bundle with an equivariant connection, then there will be a 
natural representation of $G$ on $\VV\!$, and $T^V$ will be 
$G$-invariant. See also \cite{haw1}.

If $V$ has a given inner product and a compatible connection, then 
$\VV$ will have a natural inner product, corresponding to the inner 
product of sections integrated against the volume form $\omega^n/n!$.

\section{Growth of Modules}
Since (for any $N$) the algebra $\A_N=\End\HN$ is a full matrix 
algebra, its modules are classified (modulo isomorphism) by the 
positive integers. To be precise, any $\A_N$-module can be written in 
the form $\Hom(E,\HN)$, where $\A_N$ acts only on $\HN$, 
and $E$ may be any finite-dimensional vector space; the integer 
corresponding to this module is:
\begin{definition}
$\rk[\Hom(E,\HN)] := \dim E$.
\end{definition}
This also gives a natural isomorphism $\rk : K_0(\A_N) 
\isom \Z$.

 \begin{thm}\label{growth}
Let $\VV$ be any quantization of a vector bundle $V\!$, and $V_N$ the 
restriction of $\VV$ to an $\A_N$-module.  For all sufficiently large 
values of $N$,
\beq
\rk V_N = \intM \ch V \wedge \td \TM \wedge e^{N\omega/2\pi-c_1(L_0)}
\mbox.\label{rank.formula}\eeq
\end{thm}
\begin{proof}
The uniqueness result of Thm.~\ref{uniqueness} implies that any 
analytic formula for $\rk V_N$ for large $N$ must apply to any 
quantization of $V\!$.  It is therefore sufficient to look at the 
specific quantization constructed in Sec.~\ref{VV.def}.

By Lemma~\ref{VV.preliminary}, for $N$ sufficiently large, $V_N = 
\tilde V_N = \Hom(\tilde\lef^V,\HN)$; hence, $\rk V_N = \dim 
\tilde\lef^V$.  This is the kernel of the $V^*\otimes\LN$-twisted 
Dolbeault operator, and has the same dimension as the kernel of the 
$V\otimes\bar\LN$-twisted anti-Dolbeault operator.  By 
Corollary~\ref{vanish2}, this is entirely of even degree (again, for 
$N$ sufficiently large); hence, $\rk V_N$ is the index of the 
this anti-Dolbeault operator.  Equation \eqref{rank.formula} 
then follows from the Riemann-Roch-Atiyah-Singer theorem if we note that $\ch 
\bar\LN=e^{-c_1(\LN)}=e^{N\omega/2\pi-c_1(L_0)}$.
\end{proof}

This gives some interesting qualitative results.  Again writing 
$n:=\dim_\co \M$, the right hand side of Eq.~\eqref{rank.formula} is a 
polynomial in $N$ of degree $n$.  The coefficients of this polynomial 
give $n+1$ components of the Chern character of $V\!$.  The growth of 
$\rk V_N$ thus gives some --- but not in general all --- topological 
information about the bundle $V\!$.  Evidently, the sequence of 
modules $V_N$ does not carry all the information of $\VV$; there is 
important information contained in the way these modules fit together 
as $N\to\infty$.

Since $\AA$ is a quantization of the trivial line bundle, 
Eq.~\eqref{rank.formula} implies the formula 
\beq
\dim \HN = \intM \td \TM \wedge e^{N\omega/2\pi-c_1(L_0)}
\mbox,\label{dim.formula}\eeq
which, thanks to the Kodaira vanishing theorem, holds for all $N>0$.  
Comparing \eqref{rank.formula} with \eqref{dim.formula} shows that 
$\rk V_N \approx \rk V \cdot \dim \HN$, with corrections of order 
$N^{n-1}$ ($\rk V$ is the fiber dimension).

A trivial vector bundle over $\M$ can be quantized to a free module of 
$\AA$.  In that case, $\rk V_N$ must be an integer multiple of $\dim 
\HN$, but in general the deviation from this reflects the 
nontriviality of a vector bundle.

It is especially interesting to quantize a spinor bundle.  
Since $\M$ is symplectic, it is even dimensional, and spinors decompose 
as $S=S^+\oplus S^-$ into left and right handed parts.  The Dirac 
operator is odd; that is, it maps left spinors to right spinors and 
\emph{vice-versa}.  A ``quantized'' Dirac operator should act on the 
quantized spinor bundle, i.~e., $D_N:S_N\to S_N$.  If oddness of the 
Dirac operator is preserved, and if $S_N^+$ and $S_N^-$ are of 
different size, then the quantized Dirac operator will necessarily 
have a kernel.

Typically, $S^+_N$ and $S^-_N$ \emph{are} different.  In fact $\rk S^+_N - 
\rk S^-_N$ is independent of $N$ and equal to the Euler characteristic 
$\chi(\M)$.  The dimension of an $\A_N$-module is equal to its rank 
times $\dim\HN$, so the dimension of the kernel of a quantized Dirac 
operator for a manifold of nonzero Euler characteristic must be at 
least
\[
\dim\ker D_N\geq\left|\chi(\M)\right|\dim\HN 
\mbox.\] 
This may have dire consequences for the existence of quantized Dirac 
operators. I hope to discuss this further in a future paper.

Theorem \ref{growth} can also be expressed in terms of idempotents.
 \begin{cor}
Let $e\in M_m[\C(\M)]$ and $\tilde e\in M_m(\AA)$ be idempotents such 
that $\Po(\tilde e) = e$. For $N$ sufficiently large,
\beq
\tr \tilde e_N = \intM \ch e \wedge \td \TM \wedge 
e^{N\omega/2\pi-c_1(L_0)}
\mbox.\label{GQiform}\eeq
Here $\ch e$ is the Chern character of the bundle determined by $e$, 
and $\tilde e_N$ is the evaluation of $\tilde e$ at $N$.
\end{cor}
\begin{proof}
The idempotent $e$ defines a vector bundle $V\!$.  The module 
$\VV:=\AA^m \tilde e$ is a quantization of $V\!$.  We have $V_N = 
\A_N^m \tilde e_N$.  This gives $\rk V_N = \tr \tilde e_N$.
\end{proof}
I explore an implication of Eq.~\eqref{GQiform} in another paper 
\cite{haw2}.

\begin{appendix}
\section{Spectral Inequalities}
The line bundles $\LN$ continue to be as defined in 
Sec.~\ref{Toplitz}.  Specifically, $L_N=L^{\otimes N}\otimes L_0$, and 
$L_1$ is assumed to be positive.
 \begin{lem}
\label{spectrum.lemma}
If $V$ has a compatible connection and inner product, then 
the $V^*\otimes\LN$-twisted Dolbeault operator, $D_V$, is 
(essentially) self-adjoint and 
there exists a constant, $C$, such that 
\beq
\spec D_V^2 \subset \{0\}\cup[N-C,\infty)
\mbox.\label{spectrum1}\eeq
Moreover, for the trivial bundle $V=\co\times\M$, we can take $C<1$.
\end{lem}
 \begin{proof}
Let Latin indices denote holomorphic and barred Latin indices 
antiholomorphic directions in the tangent bundle. Using the \Kahler\ 
identity $\omega_{i\jb}=ig_{i\jb}$,
the Weitzenbock formula in this case takes the form
\beq
D_V^2 = -g^{i\jb} \nabla_i \nabla_\jb + N\delta + \hat{K}
\label{Dsquared}\mbox,\eeq
where $\delta$ is the grading operator on 
$\Omega^{0,*}(\M,V^*\otimes\LN)$,
\[
\hat{K} = i \gamma^\ib \gamma^j K_{\ib j} + \tfrac{i}2 
\gamma^i\gamma^j 
K_{ij} + \tfrac{i}2 \gamma^\ib\gamma^\jb K_{\ib\jb}
\mbox,\]
and $K$ is the curvature of $V^*\otimes L_0$.

The operator $D_V^2$ always preserves the $\Z_2$-grading of 
$\Omega^{0,*}(\M,V^*\otimes\LN)$ into even and odd parts, although 
it may not respect the full $\Z$-grading.  With respect to the 
$\Z_2$-grading, $D_V$ decomposes into $D_+ + D_-$, where $D_+$ maps 
even to odd and $D_-$ maps odd to even.

The first term of \eqref{Dsquared} is a positive operator, and 
$\delta\geq1$ when restricted to the odd subspace; therefore,
\beq
D_+D_- \geq N - C
\mbox,\label{DpDm}\eeq
where $C = \norm{\hat K}$ is sufficient.  This proves that any 
eigenvalue of $D_+D_-$ (the spectrum consists entirely of eigenvalues) 
is greater than $N - C$.

Let $\psi$ be an eigenvector of $D_-D_+$ with eigenvalue 
$\lambda\neq0$.  This implies that $D_+\psi\neq 0$.  Now, 
$D_+D_-(D_+\psi) = D_+\lambda\psi = \lambda (D_+\psi)$, so $\lambda$ 
is an eigenvalue of $D_+D_-$.  Therefore, $\lambda\geq N - C$.

For $V$ trivial, The assumption that $L_1$ is positive implies that 
$\delta+\hat K$ is strictly positive.  This means that 
$D^2>(N-1)\delta$, and so we can take $C<1$ in \eqref{DpDm}.
\end{proof}

 \begin{cor}
\label{vanish2}
Let $V$ be an arbitrary vector bundle with a connection, and $D_V$ the 
$V^*\otimes\LN$-twisted Dolbeault operator.  For $N$ sufficiently 
large, $\ker D_V$ is entirely of even degree, and for any nonzero 
$\psi\in\ker D_V$, the degree $0$ component of $\psi$ is 
nonvanishing. Identical results hold for the 
$V\otimes\bar\LN$-twisted anti-Dolbeault operator.
\end{cor}
 \begin{proof}
Assign an arbitrary inner product to $V$.  The given connection on $V$ 
can be decomposed into a connection compatible with the inner product 
and a self-adjoint potential.  Correspondingly, the Dolbeault operator 
decomposes as $D_V=D_0+iB$ where $D_0$ is a self-adjoint Dolbeault 
operator and $B$ is a self-adjoint and bounded Dirac matrix.  Using 
Eq.~\eqref{Dsquared} again gives
\[
\Re D_V^2 = D_0^2 - B^2 \geq N\delta - C - \norm{B}^2
\mbox.\]
Now assume that $N > C + \norm{B}^2$.  If $\psi$ is of strictly 
positive degree (i.~e.\ $\psi_0=0$) then $D_V^2 \psi \neq 0$, which 
implies $\psi\not\in\ker D_V$.

Because $D_V$ respects the $\mathbb{Z}_2$-grading, $\ker D_V$ must be 
the sum of even and odd parts.  However, if $\psi$ is of strictly odd 
degree, then it is of strictly positive degree.  Hence, $\ker D_V$ can 
have no odd part.
\end{proof}

Note that if $V$ is trivial, then the first statement can be 
strengthened to the classical Kodaira vanishing theorem, namely the 
fact that $\ker D$ is entirely of degree $0$ and thus is simply 
$\Gh(\M,\LN)$ --- a fact which was used in Sec.~\ref{convergence}.

Recall that in the proof of Lemma~\ref{almost.isomorphic}, $\Pi_N^V$ 
was defined as the idempotent such that $\Im \Pi_N^V = \ker D_V$ and 
$\ker\Pi_N^V = \Im D_V$.
 \begin{lem}
\label{Pi.norms}
If $V$ and $W$ are the same bundle, but with different connections, 
then 
\beq
\lim_{N\to\infty} \norm{\Pi_N^V-\Pi_N^W} = 0
\mbox.\label{norm.Pi.diff}\eeq
\end{lem}
\begin{proof}
Suppose initially that the $W$ connection is compatible with the inner 
product. This means that the associated Dolbeault operator $D_W$ and 
idempotent $\Pi_N^W$ will be self-adjoint.

Since different connections on the same bundle only differ by a 
potential, the difference $A:=D_V-D_W$ of the Dolbeault operators is 
bounded.

The idempotent $\Pi_N^V$ can be expressed in terms of $D_V$ as
\[
\Pi_N^V = \frac{1}{2\pi i} \oint_\C (z - D_V)^{-1} \,dz
\mbox,\]
where the contour of integration $\C$ encloses $0$ but no other 
eigenvalue of $D_V$.  An identical formula holds for $\Pi_N^W$ in 
terms of $D_W$.  The difference of these expressions gives
\beq
\Pi_N^V - \Pi_N^W = \frac{1}{2\pi i} \oint_\C 
(z-D_V)^{-1}A\,(z-D_W)^{-1} \,dz
\mbox.\label{Pi.diff}\eeq

Expanding $(z-D_V)^{-1}=(z-D_W-A)^{-1}$ as a power series in $A$ and 
taking the norm gives
\[
\Norm{(z-D_V)^{-1}} \leq 
\left[\Norm{(z-D_W)^{-1}}^{-1}-\norm{A}\right]^{-1}
\mbox.\]
Since $D_W$ is self-adjoint, the norm of $(z-D_W)^{-1}$ is just the 
reciprocal of the distance from $z$ to $\spec D_W$.  Equation 
\eqref{spectrum1} implies that (for some $C$)
\[
\spec{D_W} \subset 
\bigl(-\infty,-\sqrt{N-C}\bigr]\cup\{0\}\cup\bigl[\sqrt{N-C},\infty\bigr)
\mbox.\]
If we let the contour $\C$ be the circle about $0$ of radius  
$\tfrac12 \sqrt{N-C}$ (which is a good contour if $N>C-4\norm{A}^2$), 
then for $z\in\C$, $\norm{(z-D_W)^{-1}}\leq 
2(N-C)^{-1/2}$. Taking the norm of \eqref{Pi.diff} now gives
\[
\norm{\Pi_N^V-\Pi_N^W} \leq  2\,\norm{A}\, (N-C)^{-1/2} 
\left[\tfrac12 
(N-C)^{1/2} - \norm{A}\right]^{-1}
\mbox.\]
This clearly goes to $0$ as $N\to\infty$, thus proving the claim 
in this special case.

Idempotents constructed from two connections incompatible with the 
inner \fchoice{\linebreak}{} product can both be compared with one 
constructed from a connection that is compatible with the inner 
product; thus, this special case implies the more general result.
\end{proof} 
\end{appendix}

\subsection*{Acknowledgements}
I wish to thank Nigel Higson for extensive discussions, as well as 
Ranee Brylinski, Jean-Luc Brylinski, and Nicolaas Landsman.  This 
material is based upon work supported in part under a National Science 
Foundation Graduate Fellowship.  Also supported in part by NSF grant 
PHY95-14240 and by the Eberly Research Fund of the Pennsylvania State 
University.

\end{document}